\documentclass[francais]{cedram-aif}
\usepackage[english,frenchb]{babel}
\usepackage[latin1]{inputenc}
 \usepackage[T1]{fontenc}
 \usepackage{lmodern}

\newcommand{\T}{\mathbb{T}}
\newcommand{\Z}{\mathbb{Z}}

\newcommand{\ba}{\backslash}
\newcommand{\rg}{\rightarrow}

\newcommand{\Ff}{\mathcal{F}}
\newcommand{\Mm}{\mathcal{M}}

\title
[Variantes sur un théorème]
{Variantes sur un théorème de Candès,\\
Romberg et Tao}

\alttitle{Variations on a theorem of Candès, Romberg and Tao}

\author{\firstname{Jean--Pierre}  \lastname{Kahane}}

\address{Laboratoire de Math\'ematiques\\ 
Universit\'e Paris--Sud \\
91405 Orsay cedex (France)}

\email{Jean-Pierre.Kahane@math.u-psud.fr}

\begin{document}

\maketitle

\subsection*{Introduction, théorème CRT et programme}

On constate expérimentalement que certains signaux et certaines images sont reproduites de façon parfaite en n'utilisant qu'une partie de leurs transformée de Fourier, grâce au procédé que voici. Le signal, ou l'image, est représenté par une fonction $x$ à valeurs complexes définie sur un groupe discret commutatif $G$ (usuellement, $\Z_N$ ou $\Z_N^2$),  et on suppose $x\in \ell^1(G)$ (c'est automatique si $G$ est fini). Sa transformée de Fourier $\hat{x}$ est définie sur le groupe dual $\hat{G}$ (de nouveau, $\Z_N$ ou $\Z_N^2$ dans les cas usuels), et $\hat{x}$ appartient à $A(\hat{G})$, l'algèbre de Wiener de $\hat{G}$ $(A(\hat{G}) = \Ff \ell^1(G))$. Ainsi
$$
\|\hat{x}\|_A = \|x\|_{\ell^1}\,.
$$
Soit $\Omega$ une partie de $\hat{G}$ (la notation traduit le fait que, pour un signal, $\Omega$ est un ensemble de fréquences $\omega$). L'ensemble des $\hat{y} \in A(\hat{G})$ qui sont égaux à $\hat{x}$ sur~$\Omega$
$$\{ \hat{y} : \hat{y}|_\Omega = \hat{x}|_\Omega\}
$$
est un convexe fermé dans $A(\hat{G})$. On cherche si ce convexe admet un point unique de norme minimum et si ce point est $\hat{x}$ ; si c'est le cas, on a reconstruit $x$ par un procédé d'analyse convexe, qui est praticable. C'est le procédé réellement utilisé.

Ce procédé marche si et seulement si

\vskip2mm

\noindent $(i)$ $\hat{x}$ \textit{est le prolongement minimal de $\hat{x}|_\Omega$ dans $A(\hat{G})$}

\vspace{2mm}

Précisons que \og\textbf{le} prolongement minimal\fg\  signifie non seulement que $\hat{x}$ est \textbf{un} prolognement minimal, mais que ce prolongement minimal est \textbf{unique}.

\noindent Le présent article est consacré à l'exploration de cette condition.

\vspace{2mm}
La première contribution majeure dans cette direction est due à Candès, Romberg et Tao (théorème~1.3 de \cite{CRT}, théorème~2.1 de \cite{C}). Dans leur modèle, les temps, $t$, sont des entiers modulo $N$, de même que les fréquences, $\omega$.~Ainsi
$$
t\in G= \Z_N = \Z/N\Z\,, \qquad \omega \in \hat{G}
$$
et $\hat{G}$ est une autre copie de $\Z_N$. La dualité entre les groupes cycliques $G$ et $\hat{G}$ s'exprime~par
$$
<\omega,t>=e\Big(\frac{\omega t}{N}\Big)\,, \qquad e(u)=\exp (2\pi iu)\,.
$$
Le signal est une fonction $x(t)$ $(t\in G)$ portée par une partie $S$ de $G$ de cardinal~$s$ :
$$
t\in G\ba S \Longrightarrow x(t)=0\,, \qquad |S| =s \,.
\leqno(1)
$$
Sa transformée de Fourier, $\hat{x}(\omega)$ $(\omega\in \hat{G})$, s'écrit
$$
\hat{x}(\omega) = \frac{1}{\sqrt{N}} \sum_{t\in G} x(t) e \Big(\frac{-t\omega}{N}\Big)
$$
et l'on a
$$
{x}(t) = \frac{1}{\sqrt{N}} \sum_{\omega\in \hat{G}} \hat{x}(\omega) e \Big(\frac{t\omega}{N}\Big)\,.
$$
Dans le cas général, on a besoin des $N$ valeurs des $\hat{x}(\omega)$ pour reconstituer le signal $x(t)$ $(t\in G)$.

Le théorème de Candès, Romberg et Tao montre comment procéder à partir d'un ensemble de fréquences bien plus réduit quand on fait l'hypothèse (1). C'est le paradigme de l'échantillonnage parcimonieux (\textit{compressive sampling}, ou \textit{compressed sensing}).

\begin{theo*}[CRT] 
On suppose que le signal $x$ est porté par $s$ points (formule $(1)$). Soit $\Omega$ la partie aléatoire de $\hat{G}$ qu'on obtient en disposant au hasard suivant la probabilité naturelle $|\Omega|$ points $\omega$ de $\hat{G}$, et supposons
$$
|\Omega| > Cs \log N \qquad \mathit{avec}\qquad C=22(M+1)\,.
$$
Alors $(i)$ a lieu avec une probabilité $1-O(N^{-M})$ $(N\rg \infty)$ et la reconstitution de $x$ se fait suivant le procédé indiqué :
$$
\|x\|_1 = \inf \{\|y\|_1 : \hat{y}|_\Omega = \hat{x}|_\Omega\}
$$
et cette égalité caractérise $x[C]$. 
\end{theo*}

La démonstration donnée par les auteurs est assez difficile. Celle que je vais donner me paraît bien plus facile, et je l'introduirai en répondant d'abord aux questions que voici.

Dans le théorème CRT, il s'agit de reconstituer \textbf{un} signal $x$ à partir de $\hat{x}|_\Omega$. Peut--on faire de même pour \textbf{tous} les signaux portés par un même support $S$ ? Mieux, pour \textbf{tous} les signaux portés par $s$ points ? Comment choisir $\Omega$, et quelles sont les estimations de probabilités correspondantes ? Ce sera l'objet des variantes.

Telles sont les questions que je m'étais posées après une lecture trop partielle de \cite{C}. Depuis la rédaction de cet article, Albert Cohen m'a informé que ces questions avaient été traitées dès 2006  par Candès et Tao : avec une probabilité voisine de 1, la reconstruction est possible pour tous les signaux $x$ portés par $s$ points moyennant une condition du type
$$
|\Omega| > C s\, \log^6 N\,.
$$
Le referee a précisé les principaux résultats avec leurs références, en particulier \cite{C,CRT,CRT2,CT}. Déjà \cite{C} signale que l'on peut remplacer l'exposant 6 par 4. Le terme utilisé est \og universal encoding strategies\fg\ et le concept de \og restricted isometries\fg\ introduit dans ce but a donné lieu à de nombreux travaux. Le résultat que j'obtiens sur cette question, avec une condition du type
$$
|\Omega| > C s^2\, \log N\,.
$$
est donc peu intéressant. La méthode que je vais présenter est très différente de celles des auteurs cités et peut intéresser certains lecteurs. Elle m'amène d'ailleurs à retrouver le théorème de Candès, Romberg et Tao avec la condition
$$
|\Omega| > C s\, \log N\,.
$$

J'indiquerai d'abord une suite de propositions qui impliquent $(i)$ :
$$
(i) \Longleftarrow (ii) \Longleftarrow (iii) \Longleftarrow (iv)
$$
Les propositions $(ii)$ et $(iii)$ ne mettent en jeu que $S$ et $\Omega$, la proposition $(iv)$ seulement $s$ et $\Omega$. La proposition $(iv)$ donnera tout de suite une réponse déterministe (sans aléa) à la question de reconstruction de tous les signaux $x$ portés par $s$ points par prolongement minimal de $\hat{x}|_\Omega$  dans $A(\hat{G})$. Son exploitation sera faite suivant la méthode aléatoire utilisée en \cite{CRT}. Puis cette méthode aléatoire sera appliquée pour obtenir $(ii)$ et $(iii)$ avec une probabilité explicite, et retrouver le théorème CRT comme conséquence de la dernière variante, $V''3$. La variante déterministe est $V1$, les variantes mettant en jeu $s$ et un $\Omega$ aléatoire sont $V2$, $V'2$ et $V''2$, les variantes concernant $S$ et $\Omega$ aléatoire sont $V3$, $V'3$ et $V''3$.

La variante $V4$ indique une application de $(i)$ dans une autre direction : $G = \Z$, $\hat{G}=\T$, $\Omega=I$, un sous--intervalle de $\T$, et on suppose que le spectre $\{\lambda_n\}_{n\in \Z}$ de $\hat{x}$ est assez dispersé, au sens que $\lambda_{n+1} -\lambda_n \ge d\ge 10$. Si $|I|$ est convenable, on peut reconstruire $\hat{x}$ à partir de $\hat{x}|_I$ grâce à~$(i)$.

\subsection*{Enoncés déterministes, variante $V1$}

Il sera commode d'écrire $\int$ pour la sommation sur $\Z_N$. Ainsi
$$
\int x = \sum_{t\in G} x(t)\,.
$$
Rappelons que
$$
\|\hat{x}\|_A = \|x\|_1 = \int |x|\,.
$$
J'écrirai $\ell(G)$ (c'est aussi $\ell'(G)$ et $\ell^\infty(G)$) l'ensemble des fonctions à valeurs complexes définies sur $G$. Désormais

\vspace{2mm}

\noindent$(i)$ $\hat{x}$ \textit{est le prolongement minimal de $\hat{x}|_\Omega$ dans $A(G)$}

\vspace{2mm}

\noindent signifie à la fois une propriété de $x\in \ell(G)$ et de $\Omega \subset \hat{G}$ et le procédé correspondant pour la reconstruction de $x$ à partir de $\hat{x}|_\Omega$. Nous supposerons toujours que $x$ est porté par $S$ et que $|S|=s$ (formule~(1)).

Soit $z\in \ell(G)$ telle que $\hat{z}|_\Omega =0$, et $p \in \ell(G)$ telle que $\hat{p}|_{\hat{G}\ba \Omega} = 0$. Donc
$$
\int_G p \bar{z} = \int_{\hat{G}} \hat{p} \bar{\hat{z}} = 0\,.
$$
$(i)$ signifie que, pour tout $z\not= 0$,
$$
\int_S |x+z| + \int_{G\ba S} |z| > \int_S |x|\,.
$$
La méthode de CRT consiste à prendre
$$
\left\{
\begin{array}{lll}
p\bar{x} = |x| &\mathrm{sur} &S\\
\noalign{\vskip2mm}
|p| <1 &\mathrm{sur} &G\ba S
\end{array}
\right. \leqno(2)
$$
Alors, tenant compte de $z\not=0$,
$$
\int_S |x+z| + \int_{G\ba S}|z| > |\int_S p(\overline{x+z}) + \int_{G\ba S} p \bar{z}| = |\int_S p \bar{x}| = \int_S |x|
$$
donc $(i)$ a lieu.

Tout le travail de CRT, qui met en jeu les matrices aléatoires, consiste à construire une fonction $p$ aléatoire vérifiant (2) avec une probabilité voisine à 1. Nous allons procéder de manière différente.

Comme
$$
\int_S |x+z| \ge \int_S |x| - \int_S |z|\,,
$$
la condition

\vspace{2mm}

\noindent $(ii)$ \textit{pour toute $z\in \ell(G) \not= 0$, telle que $\hat{z}|_\Omega=0$,}
$$
\int_{G\ba S} |z| > \int_S|z|\,,
$$
entraîne $(i)$ pour tous les $x$ portés par $S$.

Soit $\lambda\in \ell(S)$ une fonction de module 1 $(|\lambda(t)|=1)$ telle~que
$$
\lambda \bar{z} = |z| \qquad \mathrm{sur}\ \ S\,.
$$
Au lieu de (2), supposons que pour un nombre $\alpha \in ]0,1[$
$$
\left\{
\begin{array}{lll}
|p-\lambda| < \alpha &\mathrm{sur} &S\\
\noalign{\vskip2mm}
|p| < 1- \alpha &\mathrm{sur} & G\ba S\,.
\end{array}
\right.
\leqno(3)
$$
Alors
$$
|\displaystyle \int_S (p \bar{z} - \lambda \bar{z})| < \alpha \int_S |\bar{z}|$$
en supposant le second membre $\not=0$, d'où en tous cas
$$\displaystyle|\int_S p \bar{z}| > (1-\alpha) \int_S |\bar{z}|$$
$$|\displaystyle \int_{G\ba S} p\bar{z} | \leq (1-\alpha) \int_{G\ba S} |\bar{z}|\,,
$$
et, compte tenu de
$$
\int_S + \int_{G\ba S} p \bar{z} = 0\,,
$$
on obtient
$$
\int_S |\bar{z}| < \int_{G\ba S}|\bar{z}|\,.
$$
Donc, pour avoir la propriété $(ii)$, il suffit que

\vspace{2mm}

\noindent$(iii)$ \textit{pour toute $\lambda \in \ell(S)$ de module $1$, il existe $p\in \ell(G)$ tel que $\hat{p}$ soit porté par $\Omega$ et qu'on ait $(3)$ pour un $\alpha\in ]0,1[$.}

\vspace{2mm}

Le choix de $\alpha = \frac{1}{2}$ va nous convenir pour le moment. Introduisons l'idempotent
$$
K(t) = \sum_{\omega\in\Omega} e\Big(\frac{\omega t}{N}\Big) \leqno(4)
$$
et cherchons $p$ sous la forme
$$
p(t) = \sum_{t'\in S} \lambda(t') \frac{K(t-t')}{K(0)}\,.
\leqno(5)
$$
$\hat{p}$ est bien porté par $\Omega$, et on a
$$
\left\{
\begin{array}{lll}
|p(t) -\lambda(t)| \le \displaystyle \sum_{t'\in S, t'\not= t} \frac{|K(t-t')|}{K(0)} &\mathrm{si} &t\in S\\
\noalign{\vskip2mm}
|p(t)| \le \displaystyle \sum_{t'\in S} \frac{|K(t-t')|}{K(0)} &\mathrm{si} &t\in G\ba S
\end{array}
\right.
\leqno(6)
$$
donc (3) est vérifié avec $\alpha = \frac{1}{2}$ sous la condition que

\vspace{2mm}

$(iv)$ \textit{pour tout $t\not=0$ on a}
$$
|K(t)| < \frac{1}{2s} K(0)\,.
\leqno(7)
$$

Observons que $(i)$ est une condition sur $x$ et $\Omega$, $(ii)$ et $(iii)$ des conditions sur $S$ et $\Omega$, qui entraînent $(i)$ pour tous les $x$ portés par $S$, et $(iv)$ une condition sur $s$ et $\Omega$, qui entraîne $(i)$ pour tous les  $x$ dont le support a $s$ points au plus. On est ainsi parvenu, très facilement, à une réponse déterministe à la dernière question posée.

\vspace{2mm}

$V1$\quad \textit{Si $\Omega$ a la propriété que l'idempotent $K$ défini en $(4)$ vérifie la condition $(iv)$, la propriété $(i)$ est valable pour toutes les fonctions $x$ portées par un ensemble de $s$ points.}

\vspace{2mm}

La construction d'idempotents $K$ vérifiant $(iv)$, quand $N$ et $s$ sont donnés, n'est pas immédiate. Elle impose que $|\Omega|$ soit assez grand. En effet,
$$
\begin{array}{ll}
N|\Omega| &= \displaystyle\sum_{\omega\in \hat{G}} |\hat{K}(\omega)|^2 = \sum_{t\in G} |K(t)|^2\\
\noalign{\vskip2mm}
&\le |\Omega|^2 + \dfrac{|\Omega|^2}{4s^2}(N-1)
\end{array}
$$
donc
$$
|\Omega| \ge 4s^2 \frac{N}{N+4s^2-1}\,.
\leqno(8)
$$
Par exemple, si $s=2$ et $N\ge 100$, $|\Omega| \ge 16$.
L'hypothèse que $K$ est idempotent n'est d'ailleurs pas essentielle, il suffit que $\hat{K}$    soit porté par  $\Omega$. Mais elle intervient de façon naturelle dans la suite.

\subsection*{Sélection aléatoire de $\Omega$, variante $V2$}

L'étape suivante est empruntée à CRT. Elle consiste à considérer $\Omega$ comme aléatoire, à savoir comme le résultat d'une sélection aléatoire sur $\Z_N$ de paramètre $\tau$. Cela veut dire que les variables aléatoires $X_n = 1_{n\in \Omega}$ sont indépendantes,~et
$$
P(X_n=1)=\tau\,, \qquad P(X_n=0)=1-\tau\,.
$$
Alors $K$ s'écrit
$$
K(t) = \sum_{n\in \Z_N} X_n\ e\Big(\frac{nt}{N}\Big)\,.
$$

Pour chaque $t\not= 0$, soit $P(t)$ la probabilité pour que (7) n'ait pas lieu :
$$
P(t) =P (|K(t)| \ge \frac{1}{2s}K(0))\,.
$$
Choisissons un entier $\nu\ge 3$. Etant donné $t$, il existe un $\varphi \in \{\frac{2k\pi}{\nu}\,,\ k=1,2,\ldots,\nu\}$ tel~que
$$
Re(K(t) e^{-i\varphi}) \ge |K(t)| \cos \frac{\pi}{\nu}\,.
$$
Posons $a=\cos \frac{\pi}{\nu}$ et
$$
P(t,\varphi) = P (Re(K(t) e^{-i\varphi}) \ge \frac{a}{2s} K(0))\,.
$$
Ainsi
$$
P(t) \le \nu \sup\limits_\varphi P(t,\varphi)\,.
$$
Posons
$$
\begin{array}{lll}
Y &= &Re(K(t)e^{-i\varphi}) - \frac{a}{2s}	K(0)\\
\noalign{\vskip2mm}
&=&\displaystyle \sum_{n\in \Z_N} X_n \Big(\cos\Big(\frac{2\pi nt}{N} - \varphi\Big) -\frac{a}{2s}\Big)\\
\noalign{\vskip2mm}
&=&\displaystyle \sum_{n\in \Z_N} X_n A_n\,.
\end{array}
$$
On a
$$
P(t,\varphi) = P(Y\ge 0) < E(e^{uY})
$$
pour tout $u>0$. Majorons $E(e^{uY})$. Comme
$$
E(e^{vX_n}) = 1+\tau (e^v-1) < \exp (\tau(e^v-1))
$$
et
$$
E(e^{u\Sigma A_nX_n}) = \Pi E(e^{uA_nX_n})\,,
$$
on a
$$
E(e^{uY}) \le \exp (\tau \sum_{n\in \Z_N} (\exp (uA_n)-1))\,.
$$
Désignons par $\Mm$ la moyenne sur $\Z_N : \Mm = \frac{1}{N} \sum\limits_{n\in\Z_N}$, et posons
$$
\begin{array}{lll}
B_n &=& \cos \Big(\dfrac{2\pi nt}{N}-\varphi\Big)\,,\\
\noalign{\vskip2mm}
A_n&= &B_n -\dfrac{a}{2s}\,.
\end{array}
$$
Il vient
$$
E(e^{uY}) \le \exp (\tau N(\exp\Big(\frac{-au}{2s}\Big) \Mm(e^{uB_n})-1))\,.
$$
Ecrivons
$$
\Mm(e^{uB_n}) = \sum_{k=0}^\infty \frac{u^k}{k!} \Mm(B_n^k)
$$
et évaluons $\Mm(B_n^k)$.

Comme $t\not=0$,
$$
\Mm(B_n)=0\,.
$$
Si de plus $2t\not= 0$,
$$
\Mm(B_n^2) = \frac{1}{2}\,.
$$
Si de plus $3t\not=0$,
$$
\Mm(B_n^3)=0\,.
$$
Supposons qu'il en est bien ainsi, c'est--à--dire que $N$ n'est multiple ni de 2 ni de 3. En tous cas
$$
\Mm(B_n^k)\le 1\,.
$$
Donc
$$
\Mm(e^{uB_n}) \le 1+ \frac{u^2}{4} + \sum_{k=4}^\infty \frac{u^k}{k!}
$$
et
$$
\begin{array}{rl}
E(e^{uY}) &\displaystyle\le \exp\Big(\tau N\Big(\Big(\Big(1- \dfrac{au}{2s} + \dfrac{a^2u^2}{8s^2}\Big)\wedge 1\Big) \Big(1+ \dfrac{u^2}{4} + \sum_{k=4}^\infty \dfrac{u^k}{k!}\Big) -1\Big)\Big)\\
&\le \displaystyle\exp\Big(\tau N\Big(- \dfrac{au}{2s} + \dfrac{a^2u^2}{8s^2} + \dfrac{u^2}{4} + \sum_{k=4}^\infty \dfrac{u^k}{k!}\Big)\Big)\,.
\end{array}
$$
Le choix de $u=\frac{a}{s}\Big(1+ \frac{a^2}{s^2}\Big)^{-1}$ donne
$$
-\frac{au}{2s} + \frac{a^2u^2}{8s^2} + \frac{u^2}{4} + \sum_{k=4}^\infty \frac{u^k}{k!} < - \frac{au}{2s} + \frac{u^2}{4} \Big(1+ \frac{a^2}{s^2}\Big) = - \frac{a^2}{4s^2}\Big(1+ \frac{a^2}{s^2}\Big)^{-1}
$$
donc
$$
\begin{array}{l}
P(t,\varphi) = P(Y\ge 0) < \exp \Big(-\tau N \dfrac{a^2}{4(s^2+a^2)}\Big)\,,\\
P(t) < \nu \exp \Big( -\tau N \dfrac{a^2}{4(s^2+a^2)}\Big)\,,\\
P(non\, iv)) < \nu N\exp \Big(-\tau N \dfrac{a^2}{4(s^2+a^2)}\Big)\,.
\end{array}
$$
Posons
$$
\tau N = 4C (s^2+1) \log N\,;
$$
alors
$$
P(non\,(iv)) < \nu N^{-Ca^2+1}\,.
$$

Par exemple, si $N=1001$, $s=2$, $\tau N=300$ et $\nu=10$, on obtient $P((iv)) > \frac{99}{100}$. 
On a donc un moyen pratique de réaliser $(iv)$ par tirage au sort, mais on doit imposer à $|\Omega|$ des valeurs assez grandes $(E(|\Omega|) = \tau N)$. Enonçons le résultat.

\vspace{2mm}

$V2$ \ \ \textit{Supposons qu'on réalise $\Omega$ par sélection aléatoire dans $\Z_N$, avec}
$$
E(|\Omega|) = 4C (s^2+1) \log N\,, \qquad C>1\,,
$$
\textit{et la condition que $N$ ne soit multiple ni de $2$, ni de $3$. Si $\delta=Ca^2-1$ avec $a=\cos \frac{\pi}{\nu}$, $\nu\ge 3$,}
$$
P((iv)) >1 - \nu N^{-\delta}\,.
$$

Rappelons que, si $(iv)$ a lieu, $(i)$ est valable pour \textbf{toutes} les  fonctions $x$ portées par un ensemble de $s$ points.

\subsection*{Deux probabilités sur les $\Omega$, variantes $V'2$ et $V''2$}

Avant d'aller plus loin, revenons sur les deux probabilités qui ont été introduites sur l'ensemble des parties de $\Z_N$. La première, qui figure dans l'énoncé du théorème CRT, est définie par la donnée d'un entier $f$ (l'effectif des fréquences), $0<f<N$ ; elle est concentrée sur l'ensemble des parties $\Omega$ de $\Z_N$ de cardinal $|\Omega|=f$ et elle y est également distribuée. La seconde intervient dans la variante $V2$, elle est définie par un (petit) paramètre $\tau$, $0<\tau<1$ ; les événements $(n\in \Omega)$ $(n\in N)$ sont indépendants et ont $\tau$ pour probabilité commune ; on dit qu'$\Omega$ résulte d'une sélection aléatoire de paramètre $\tau$. On peut également la définir par le fait que la loi des $|\Omega|$ est la distribution de Bernoulli $B(\tau,N)$ et que, pour chaque $0\le j\le N$, la distribution des $\Omega$ tels que $|\Omega|=j$ est uniforme (égale probabilité sur tous les $\Omega$). Dans les deux cas, la probabilité est invariante par les permutations de $\Z_N$, et elle est donc bien définie par la distribution des $|\Omega|$. Désignons ces probabilités par $P_{1,f}$ et~$P_{2,\tau}$ :
$$
\begin{array}{l}
P_{1,f}(|\Omega|=f)=1\\
P_{2,\tau}(|\Omega|=j)=\begin{pmatrix}
N\\j
\end{pmatrix} \tau^j (1-\tau)^{N-j}\,.
\end{array}
$$
Il est classique, et facile à voir, que
$$
E_{2,\tau}(|\Omega|) = \tau N
$$
et
$$
P_{2,\tau} (|\Omega| > \tau N(1+\varepsilon)) < \exp \Big(-\frac{\varepsilon^2}{3} \tau N \Big) \qquad \mathrm{si}\ 0<\varepsilon<\frac{1}{10}\,. \leqno(9)
$$

Les propriétés $(i)$, $(ii)$ et $(iii)$ sont monotones en $\Omega$ : si elles ont lieu avec un $\Omega$, elles ont lieu avec un $\Omega' \supset \Omega$. Cela permet de s'affranchir des conditions arithmétiques indiquées en $V2$. Il n'est pas clair qu'il en soit de même pour $(iv)$. Si l'on s'en tient à $(i)$, $(ii)$ et $(iii)$, (9) montre qu'en choisissant $f= \tau N(1+\varepsilon)$ on a
$$
\begin{array}{ll}
P_{1,f}(\cdot) &\ge P_{2,\tau} (\cdot) - \exp\Big(-\dfrac{\varepsilon^2}{3} \tau N\Big)\\
\noalign{\vskip2mm}
&> P_{2,\tau} (\cdot) - \exp\Big(-\dfrac{\varepsilon^2}{4} f\Big) \qquad \mathrm{si}\ \varepsilon < \dfrac{1}{12}\,.
\end{array}
$$
La variante $V2$, où $P(\cdot)$ signifie $P_{2,\tau}(\cdot)$, a donc pour corollaire :

\vspace{2mm}

$V'2$\ \ \textit{Si $\Omega$ est choisi au hasard parmi les parties de $\Z_N$ de cardinal $f$, avec la probabilité naturelle, et~si}
$$
f=4C(1+\varepsilon)(s^2+1) \log N\,, \qquad C>1,\ O<\varepsilon<\frac{1}{12}\,,
$$
\textit{la probabilité d'avoir $(i)$ pour tous les $x$ portés par un ensemble de $s$ points, qu'on va désigner par $P(\forall x \mid \|x\|_0 \le s$ $(i))$, vérifie}
$$P(\forall x \mid \|x\|_0 \le s \ (i)) >1 - \nu N^{-\delta} - N^{-C \varepsilon^2(1+\varepsilon)(s^2+1)}
$$
\textit{avec $\delta = Ca^2 -1$, $a=\cos \frac{\pi}{\nu}$, $\nu \ge 3$ arbitraire.}

\vspace{2mm}

Ici $P=P_{1,f}$, et $\|x\|_0$ désigne le cardinal du support de~$x$.

L'énoncé est explicite mais lourd. Pour le comparer au théorème CRT, en voici une version allégée, asymptotique :

\vspace{2mm}

$V''2$\ \ \textit{Si $\Omega$ est choisi au hasard parmi les parties de $\Z_N$ de cardinal $f$, et~si}
$$
f =|\Omega| = 4C(s^2+1) \log N\,, \quad C>1\,,
$$
\textit{la probabilité d'avoir $(i)$ pour tous les $x$ portés par un ensemble de $s$ points vérifie}
$$
P( \forall x \mid \|x_0\| \le s\ \ (i)) = 1-O(N^{-\delta}) \qquad (N\rg \infty)
$$
\textit{pour tout $\delta < C-1$.}

\vspace{2mm}

En ce qui concerne la probabilité d'avoir $(i)$ pour un $x$  fixé, cet énoncé n'entraîne  le théorème CRT que dans le cas $s=2$, avec une amélioration insignifiante. Son intérêt est de s'appliquer à tous les $x$ portés par $s$ points. Peut--on remplacer $s^2$ par une plus petite puissance de $s$ ? La question est ouverte pour $V'2$ et $V''2$. La réponse est négative pour $V2$, en vertu de l'inégalité~(8).

\subsection*{Variantes $V3$, $V'3$ et $V''3$}

Nous supposons désormais que $x$ est porté par un ensemble $S$ fixé mais inconnu, tel que $|S|=s$. Nous allons d'abord utiliser pour $\Omega$ le modèle de la sélection aléatoire de paramètre $\tau$ sur $\Z_N$, et chercher à réaliser la proposition $(iii)$, que je rappelle : il existe un $\alpha \in ]0,1[$ tel que, pour toute $\lambda \in \ell(S)$ telle que $|\lambda|=1$, il existe $p\in \ell(G)$ à spectre dans $\Omega$ vérifiant $|p-\lambda|<\alpha$ sur $S$ et $|p| <1-\alpha$ sur~$G\ba S$.

Pour chaque $t'\in S$, remplaçons $\lambda(t')$ par la plus proche racine $\mu$--\up{ième} de l'unité, $e^{i\psi(t')}$. Ainsi
$$
|e^{i\psi(t')}-\lambda(t') | \le 2 \sin \frac{\pi}{2\mu}
$$
Posons
$$
\alpha' = \alpha -2 \sin \frac{\pi}{2\mu}
$$
et supposons $\alpha'>0$. Alors
$$
|p-e^{i\psi}| < \alpha' \quad \mathrm{sur}\ \ S
$$
entraîne $|p-\lambda| <\alpha$ sur $S$. Choisissons
$$
p(t) = \sum_{t'\in S} e^{i\psi(t')} \frac{K(t-t')}{K(0)}\,,
$$
$K$ étant l'idempotent défini en (4), tel que $\hat{K}=\sqrt{N} 1_\Omega$. 
Ainsi $p$ a son spectre dans $\Omega$ et vérifie les inégalités voulues si
$$
\left\{
\begin{array}{lll}
\displaystyle \Big(\sum_{t'\in S,t'\not= t} e^{i\psi(t')} \frac{K(t-t')}{K(0)} \Big) < \alpha' &\mathrm{sur} &S\\
\displaystyle \Big(\sum_{t'\in S} e^{i\psi(t')} \frac{K(t-t')}{K(0)} \Big) <1-\alpha &\mathrm{sur} &G\ba S\,.
\end{array}
\right.\leqno(10)
$$
Soit $e^{i\varphi(t)}$ la racine $\nu$\up{ième} de l'unité la plus proche de $p(t)-e^{i\psi(t)}$ si $t\in S$, et de $p(t)$ si $t\in G \ba S$. Posons $a=\cos \frac{\pi}{\nu}$. Pour avoir (10), il suffit que
$$
\left\{
\begin{array}{lll}
Re\displaystyle \sum_{t'\in S,t'\not= t} e^{i(\psi(t')-\varphi(t))} \frac{K(t'-t)}{K(0)}  < a \alpha' &\mathrm{sur} &S\\
Re\displaystyle \sum_{t'\in S} e^{i(\psi(t')-\varphi(t))} \frac{K(t'-t)}{K(0)}  < a(1-\alpha) &\mathrm{sur} &G\ba S\,.
\end{array}
\right.
$$
c'est--à--dire
$$
\left\{
\begin{array}{ll}
\kern-2mm Y =\kern-2mm \displaystyle
\sum_{n\in \Z_N}\!\!X_n \!\!\sum_{t'\in S\ba \{t\}}\!\!\cos \big(\frac{2\pi n(t\!-\!t')}{N} + \psi(t')\! -\!\varphi(t)\big) \!-\!a\alpha' \!<\!0
 \kern-2mm &\mathrm{sur}\ S\\
\kern-2mm Z  = \kern-2mm\displaystyle\sum_{n\in \Z_N}\!\!X_n \sum_{t'\in S}\cos \big(\frac{2\pi n(t\!-\!t')}{N} + \psi(t') \!-\!\varphi(t)\!\big)-\!a(1\!-\!\alpha) \! <\!0 \kern-2mm &\mathrm{sur}\ G\ba S\,.
\end{array}
\right.
$$
Outre $N$ et $S$, fixons $\alpha$, $\mu$ et $\nu$, et évaluons la probabilité pour qu'il n'en soit pas ainsi, c'est--à--dire qu'il existe une fonction $\psi$ ($\mu^s$ choix possibles) et, soit un $t\in S$ et un $\varphi(t)$ ($s\nu$ choix possibles) tels que $Y\ge 0$, soit un $t\in G\ba S$ et un $\varphi(t)$ ($(N-s)\nu$ choix possibles) tels que $Z\ge 0$.

Evaluons $P(Y\ge 0)$ quand $\psi$, $t$ et $\varphi(t)$ sont fixés, par la même  méthode que pour obtenir~$V2$ :
$$
P(Y\ge 0) < E(e^{uY}) <\exp (\tau N(e^{-a\alpha'u}\Mm e^{uD_n}-1)) \quad (u>0)
$$
avec $\Mm =\frac{1}{N} \sum\limits_{n\in \Z_N}$ et
$$
D_n = \sum_{t'\in S \ba\{t\}} \cos \Big(\frac{2\pi n(t-t')}{N} + \psi(t') -\varphi(t)\Big)\,.
$$
On a
$$
\begin{array}{c}
e^{-a\alpha'u} < (1-a\alpha' u + \frac{u^2}{2} (a\alpha')^2) \wedge 1\\
\Mm(e^{uD_n}) = 1 +u \ \Mm(D_n) + \displaystyle \sum_{k=2}^{\infty} \frac{u^k}{k!} \Mm(D_n^k
)\,.
\end{array}
$$
Comme $t-t'\not=0$ dans $\Z_N$, on a $\Mm(D_n)=0$. Ecrivons
$$
D_n^k = \frac{1}{2^k} \sum_{\varepsilon_j=\pm1,t_j'\in S\ba \{t\},j=1,2,\ldots,k} \cos\!\!\! \sum_{j=1,2,\ldots,k}
\varepsilon_j \big(\frac{2\pi n(t-t_j')}{N} + \psi(t_j')-\varphi (t)\big).
$$
La moyenne du cosinus est nulle sauf si $\Sigma\varepsilon_j (t-t_j')=0$ dans $\Z_N$. Quand les $\varepsilon_j(j=1,2,\ldots,k)$ et les $t_j'(j=1,2,\ldots, k-1)$ sont fixés, cette égalité détermine $t_k'$. Donc
$$
\begin{array}{lll}
\Mm(D_n^k) &\le &(s-1)^{k-1} \qquad (k=2,3,\ldots)\,,\\
\Mm(e^{uD_n}) &\le &1+ \displaystyle \sum_{k=2}^\infty \frac{u^k(s-1)^{k-1}}{k!}\,.
\end{array}
$$
Imposons $u(s-1) \le 1$. Alors
$$
\begin{array}{c}
\Mm(e^{uD_n}) \le  1+u^2(s-1)(\frac{1}{2!}+\frac{1}{3!}+\cdots) = 1+u^2(s-1)(e-2)\\
e^{-a\alpha'u} \Mm(e^{uD_n})-1 \le -a\alpha' u +\frac{u^2}{2}(a\alpha')^2 +u^2 (s-1)(e-2)\,.
\end{array}\leqno(11)
$$
Le minimum est atteint pour $u((a\alpha')^2 +2(s-1)(e-2))=a\alpha'$ (d'où $u(s-1)\le 1)$ et vaut
$$
-\frac{(a\alpha')^2}{2(a\alpha')^2+4(s-1)(e-2)}
$$
donc
$$
P(Y\ge 0) \le \exp \Big(-\tau N\frac{a^2\alpha'^2}{2a^2\alpha'^2+4(s-1)(e-2)}\Big)\,.
$$
Mutatis mutandis, le même calcul donne
$$
P(Z\ge 0) \le \exp \Big( -\tau N \frac{a^2(1-\alpha)^2}{2a^2(1-\alpha)+4s(e-2)}\Big)\,,
$$
d'où finalement l'évaluation
$$
\begin{array}{c}
1-P((ii)) \le \nu\mu^s \Big( s\exp (-\tau N \dfrac{a^2\alpha'^2}{2a^2\alpha'^2+4(s-1)(e-2)}\Big) +\\
\hfill (N-s) \exp \Big( -\tau N 
\dfrac{a^2(1-\alpha)^2}{2a^2(1-\alpha)^2+4s(e-2)}\Big)\Big)\,,
\end{array}
\leqno(12)
$$
où $N$ et $s$ sont les données, et on cherche à déterminer $\tau N=E(|\Omega|)$ de façon que, par un choix judicieux des entiers $\nu$ et $\mu \ge4$ et de  $\alpha\in ]2 \sin \frac{\pi}{2\mu},1[$, avec $a=\cos \frac{\pi}{\nu}$ et $\alpha'=\alpha -2\sin \frac{\pi}{2\mu}$, on ait une majoration convenable de $1-P((ii))$.

Simplifions en remplaçant $4s(e-2)$ par $3s$ et en posant
$$
\begin{array}{c}
\tau N = (3s+2) C \log N : \\
1-P((ii)) \le \nu\mu^s(sN^{-C a^2\alpha'^2}+N^{-Ca^2(1-\alpha)^2+1})
\end{array}
\leqno(13)
$$
Simplifions encore en imposant la condition
$$
Ca^2 \alpha'^2 = Ca^2(1-\alpha)^2 -1
$$
qui s'écrit aussi
$$
Ca^2(1-2 \sin \frac{\pi}{2\mu})(1-2 \sin \frac{\pi}{2\mu} - 2\alpha')=1\,.
\leqno(14)
$$
Enonçons le résultat à ce stade.

\vspace{2mm}

$V3$\ \ \textit{Supposons qu'on réalise $\Omega$ par sélection aléatoire dans $\Z_N$, avec $E(|\Omega|) = \tau N=(3s+2)C \log N$. Choisissons $\mu$ et $\nu$ entiers $\ge 4$. Sous la condition $(14)$, la probabilité de la condition $(ii)$ vérifie}
$$
P((ii)) \ge 1- \nu \mu^s (s+1)N^{-Ca^2 \alpha'^2}\,,
$$
\textit{où $a=\cos \frac{\pi}{\nu}$.}

\vspace{2mm}

Rappelons que $(ii)$ entraîne la validité de $(i)$, c'est--à--dire la possibilité de reconstruire $x$ par extrapolation de $\hat{x}|_\Omega$ dans $A(\hat{G})$, pour tous les signaux $x$ portés par un ensemble $S$ (fixe mais non précisé), de cardinal $|S|=s$.

Pour tester $V3$, prenons $N=10^{10}$, $s=10$, $\mu=\nu=10$, et $\nu N=E(|\Omega|)=1,5\cdot 10^4$. On obtient $P((ii)) >1 -0,5 \cdot 10^{-4}$. C'est un peu meilleur que ce qu'on a en appliquant~$V2$.

L'estimation asymptotique quand $N\rg \infty$ est plus intéressante. En effet, fixons $C>1$. Pour un $\mu$ suffisamment grand, $\alpha'$ est arbitrairement proche de $\frac{1}{2}\big(1-\frac{1}{C}\big)$, donc $C a^2\alpha'^2$ de $\frac{|C-1]^2}{4C}a^2$. Pour $\nu$ suffisamment grand, $a$ est arbitrairement proche de $1$. Enonçons le résultat.

\vspace{2mm}

$V'3$\ \ \textit{Supposons qu'on réalise $\Omega$ par sélection aléatoire dans $\Z_N$, avec}
$$
E(|\Omega|) = \tau N =(3s+2)C \log N\,, \qquad C>1\,,
$$
et soit 
$$
\delta < \frac{(C-1)^2}{4C}\,.
$$
Alors
$$
P((ii)) = 1 -O (N^{-\delta}) \qquad(N\rg \infty)
$$

La comparaison que nous avons faite des deux procédés aléatoires pour réaliser $\Omega$ donne en corollaire :

\vspace{2mm}

$V''3$\ \ \textit{Soit $\Omega$ la partie aléatoire de $\hat{G}=\Z_N$ qu'on obtient en disposant au hasard suivant la probabilité naturelle $|\Omega|$ points $\omega$ de $\hat{G}$, et supposons}
$$
|\Omega| = (3s +2) C \ \log N
$$
\textit{avec $C>1$ et}
$$
\delta < \frac{(C-1)^2}{4C}\,.
$$
\textit{Alors la probabilité de pouvoir utiliser $(i)$ pour la reconstruction de tous les signaux $x$ portés par un ensemble $S$ (fixe mais inconnu) est $1-O(N^{-\delta})$ $(N\rg \infty)$.}

On a ainsi retrouvé le théorème CRT, avec une estimation des constantes, et une démonstration que j'espère facile à comprendre sinon plaisante à suivre dans le détail des calculs.

On peut d'ailleurs améliorer les constantes et affiner les calculs. Pour simplifier l'écriture, j'ai remplacé à partir de (13) $4s(e-2)+2$ par $3s+2$. Dans chacun des énoncés $V3$, $V'3$, $V''3$, on pourrait écrire $2,874\ s +2$ à la place de $3s+2$. On peut aller un peu plus loin en imposant $u(s-1)\le \xi$~d'où
$$
\Mm(e^{uD_n}) \le 1+ u^2 (s-1)\big(\frac{1}{2!}+\frac{\xi}{3!}+\frac{\xi^2}{4!}+\cdots\big)
$$
au lieu de (11), et choisir $\xi$ de façon que
$$
u(a\alpha')^2 + 2(s-1)\big(\frac{1}{2!}+\frac{\xi}{3!}+\frac{\xi^2}{4!} +\cdots\big) = a\alpha'
$$
entraîne $u(s-1) \le \xi$. Le choix de $\xi=0,763$ est convenable et permet de remplacer $2,874$ par $2,623$. Tant qu'on n'a pas de limitation dans l'autre sens, abaisser ainsi le coefficient de $s$ dans les énoncés $V3$, $V'3$ et $V''3$ est un jeu gratuit.

Sous des formes voisines, les variantes de $V1$ à $V''3$ sont énoncés dans~\cite{K'}.

\subsection*{Variante $V4$}

Il y a d'autres variations possibles autour de la propriété $(i)$. Je pense à la détermination d'une fonction $\in A(\T)$ dont le spectre, inconnu, est assez dispersé, à partir de sa donnée sur un intervalle. Voici un résultat typique de ce qu'on peut chercher dans cette direction.

\vspace{2mm}

$V4$\ \ \textit{Soit $\hat{x} \in A(\T)$ une fonction dont le spectre $\{\lambda_n\}_{n\in \Z}$ a un pas\break $\ge d\ge 10$ :}
$$
\begin{array}{l}
\hat{x}(t) = \displaystyle \sum_{n\in \Z} x(\lambda_n) e^{2\pi i\lambda_n t}\\
\displaystyle \sum_{n\in \Z} |x(\lambda_n)| < \infty\,, \quad \lambda_{n+1} - \lambda_{n} \ge d \quad (n\in \Z)\,.
\end{array}
$$
\textit{Soit $I$ un intervalle dans $\T$. Si}
$$
|I| \ge \frac{10}{d} \log d\,,
$$
\textit{$\hat{x}$ est le prolongement minimal de $\hat{x}|_I$ dans $A(\T)$.}

\vspace{2mm}

La preuve et des commentaires se trouvent en \cite{K}.

\vskip4mm

Je remercie ma première lectrice, Anne de Roton, pour ses observations pertinentes et pour les corrections qu'elle m'a fait faire.

J'ai déjà signalé les informations et critiques fructueuses apportées par Albert Cohen et par le referee, que je remercie chaleureusement.



\begin{thebibliography}{3}


\bibitem{C}
Emmanuel J. \textsc{Candès}.--- \textit{Compressive Sampling}, Proceedings of the International Congress of Mathematicians, Madrid 2006, Volume III, 1433--1452.

\bibitem{CRT}
E. J. \textsc{Candès}, J. \textsc{Romberg} and T. \textsc{Tao}.--- \textit{Robust Uncertainty Principles : Exact Signal Reconstruciton From Highly Incomplete Frequency Information,}  IEEE Transactions on Information Theory 52, 2 (2006), 489--509.

\bibitem{CRT2}
E. J. \textsc{Candès}, J. \textsc{Romberg} and T. \textsc{Tao}.--- \textit{Stable signal recovery from incomplete and inaccurate measurements,}  Communications in Pure and Applied Mathematics 59, 8 (2006), 1207--1223.

\bibitem{CT}
E. J. \textsc{Candès} and T. \textsc{Tao}.--- \textit{Near optimal signal recovery from random projections : universal encoding strategies,}  IEEE Transactions on Information Theory 52, 12 (2006), 5406--5425.

\bibitem{Co}
A. \textsc{Cohen}.--- \textit{Communication orale lors du colloque de décembre 2011 au Centre de recerca matem\'atica (CRM) de l'Université autonome de Barcelone,}  ICREA Conference on Approximation Theory and Fourier Analysis.


\bibitem{K}
J.--P. \textsc{Kahane}.--- \textit{Analyse et synthèse harmonique,} Histoires de mathématiques (Journéees X--UPS 2011), \'Ecole Polytechnique, Palaiseau, 2012, 17--53.

\bibitem{K'}
J.--P. \textsc{Kahane}.--- \textit{Idempotents et échantillonnage parcimonieux}, Comptes rendus de l'Académie des sciences de Paris. série 1, 349 (2011), 1073--1076.



\vskip4mm






\end{thebibliography}
\end{document}